\documentclass{article}
\usepackage[utf8]{inputenc}
\usepackage[T1]{fontenc}
\usepackage{amsmath,amssymb,amsthm}
\usepackage{anysize}
\usepackage{url}
\usepackage[dvipsnames]{xcolor}
\usepackage{paralist}
\usepackage{graphicx}
\usepackage{pgfplots}
\pgfplotsset{compat=1.4} 
\usepackage{textcomp}
\usepackage{tikz} 
\usepackage{tkz-euclide} 
\usetikzlibrary{shapes,shadows,patterns,trees,automata,positioning,snakes,fadings,backgrounds,intersections,calc,arrows.meta} 

\usepackage{ifthen}

\theoremstyle{plain}
\newtheorem{theorem}{Theorem}[section]    
\newtheorem{conjecture}[theorem]{Conjecture} 
\newtheorem{lemma}[theorem]{Lemma}

\newtheorem{proposition}[theorem]{Proposition}
\newtheorem*{proposition*}{Proposition}
\newtheorem{corollary}[theorem]{Corollary}
\theoremstyle{definition}  
\newtheorem{definition}[theorem]{Definition}
\newtheorem{example}[theorem]{Example}
\newtheorem*{definition*}{Definition}   
\theoremstyle{remark} 
\newtheorem{remark}[theorem]{Remark} 
\newtheorem*{problem*}{Problem} 

\newcommand{\N}{\ensuremath{\mathbb N}}
\newcommand{\Z}{\ensuremath{\mathbb Z}}
\newcommand{\R}{\ensuremath{\mathbb R}}

\newcommand{\EE}{\ensuremath{\mathrm{E}}}
\newcommand{\EEM}{\ensuremath{\mathbf{EM}}}
\newcommand{\DS}{\ensuremath{\mathrm{DS}}} 
\newcommand{\DSM}{\ensuremath{\mathbf{DSM}}} 
 
\newcommand{\fset}{\mathcal{F}}

\newcommand{\aff}{\operatorname{aff}}
\newcommand{\lin}{\operatorname{lin}}
\newcommand{\pyr}{\operatorname{pyr}}
\newcommand{\bip}{\operatorname{bip}}
\newcommand{\prism}{\operatorname{prism}}
\newcommand{\dhalfup}{\lceil d/2\rceil}
\newcommand{\dhalfdown}{\lfloor d/2\rfloor} 
\newcommand{\pluszero}{\mathbin{{\scriptstyle\boxplus}}}
\newcommand{\plusone}{\mathbin{{\scriptstyle\boxplus}'}}
\newcommand{\plustwo}{\mathbin{{\scriptstyle\boxplus}''}}
 
\usepackage[colorlinks=true,urlcolor=blue,citecolor=magenta]{hyperref}
\hypersetup
{
  pdftitle   = {additive-f-vectors},
  pdfsubject = {f-vector additive},
  pdfauthor  = {Günter M. Ziegler}
}

\title{Additive structures on $f$-vector sets of polytopes} 
\author{{\Large Günter M. Ziegler\footnote{The author was supported by the DFG Collaborative Research Center TRR~109 ``Discretization in Geometry and Dynamics.''
This material is based on work
supported by the National Science Foundation under Grant No.~DMS-1440140 while
the author was in residence at the
Mathematical Sciences Research Institute in Berkeley CA, during the Fall 2017 semester.}}\\[2pt]
Institut f\"ur Mathematik, FU Berlin\\Arnimallee 2, 14195 Berlin, Germany\\
\url{ziegler@math.fu-berlin.de}}
\date{September 13, 2017}

\begin{document}
\maketitle 

\begin{abstract} 
We show that the $f$-vector sets of $d$-polytopes have non-trivial additive structure: They span affine lattices and are embedded in monoids that we describe explicitly. Moreover, for many large subclasses, such as the simple polytopes, or the simplicial polytopes,  there are monoid structures on the set of $f$-vectors by themselves: ``addition of $f$-vectors minus the $f$-vector of the $d$-simplex'' always yields a new $f$-vector.  For general $4$-polytopes, we show that the modified addition operation does not always produce an $f$-vector, but that  the result is always close to an $f$-vector. In this sense, the set of $f$-vectors of \emph{all} $4$-polytopes forms an ``approximate affine semigroup.''  The proof relies on the fact for $d=4$ every $d$-polytope, or its dual, has a ``small facet.'' This fails for $d>4$. 	 

We also describe a two further modified addition operations on $f$-vectors that can be geometrically realized by glueing  corresponding polytopes. The second one of these may yield a semigroup structure on the $f$-vector set of all $4$-polytopes.
%
\end{abstract}

\section{Introduction}

Let $P$ be a $d$-dimensional convex polytope, or \emph{$d$-polytope} for short. 
(For introductions and surveys see \cite{Gr1-2}, \cite{Z35}, or \cite{Z49-3}.)
The \emph{vector of face numbers} $f(P)$, known as the \emph{$f$-vector} for short, is the vector
$f(P):=(f_0,\dots,f_{d-1})\in\Z^d$, where $f_i$ denotes the number of $i$-dimensional
faces of~$P$. 
For example, the $f$-vector of a $d$-simplex,
$f(\Delta_d)=\big(d+1,\binom{d+1}2,\dots,\binom{d+1}{d-1},d+1\big)$,
is a row of the Pascal triangle.

Let $\mathcal{P}^d$ denote the set of combinatorial types of $d$-dimensional polytopes
and write
\[
	\fset(\mathcal{P}^d):=\{(f_0(P),\dots,f_{d-1}(P)):P\in\mathcal{P}^d\}\subset\Z^d
\]
for the set of all $f$-vectors of $d$-dimensional polytopes.
We will write $\mathcal{P}^d_s$ for the set of combinatorial
types of simplicial $d$-polytopes and $\fset(\mathcal{P}^d_s)$
for its $f$-vector set, and analogously for other classes of $d$-polytopes.

The sets $\fset(\mathcal{P}^d)$ of all $f$-vectors are
easy to determine completely for $d=2$ and $d=3$, and very structured:
In both cases we get the set of all integer points in an affine polyhedral
cone of dimension $d-1$, and thus the structure of a 
finitely-generated affine monoid.

The affine hull of the set of $f$-vectors of all $d$-polytopes
is $(d-1)$-dimensional for all $d\ge2$.
However, for $d\ge4$ the sets $\fset(\mathcal{P}^d)$ seem complicated
(even provably so, see Sjöberg \& Ziegler \cite{Z163}),
and only very partial information is available. In particular, 
for $d\ge4$ the set $\fset(\mathcal{P}^d)$ is not the set of all 
integer points in a convex set,
as constraints like $f_1\le\binom{f_0}2$ are concave rather than convex.

Nevertheless, this paper presents the (surprising, at least to the author) 
observation that there is a natural structure of an 
(additive, commutative, cancellative, usually not finitely-generated)
affine monoid or semigroup on various $f$-vector sets $\fset$ of $d$-polytopes.

For this we define an \emph{affine semigroup} as any structure
that is isomorphic to a subset of some integer
lattice $\Z^d$ that is closed under addition. It is an
\emph{affine monoid} if it additionally has a neutral element.
Such a semigroup is \emph{cancellative} no element,
other than a possible neutral element, has an inverse.
For such a semigroup the minimal (with respect to inclusion)
set of generators is unique, but it need not be finite.
(It is also known as the \emph{Hilbert basis}.)
All the semigroups we consider are isomorphic via some integral
translation to subsemigroups of some $(\Z^d,+)$.
Hence in particular they are commutative.

The monoid structures on $f$-vector sets $\fset$ that we construct 
arise from the addition operation on $\Z^d$ defined by
\[
	x\pluszero y:=x+y-f(\Delta_d),
\]
so $f(\Delta_d)$ is a neutral element for this.
Other semigroup structures that we obtain use the addition operations defined by
\[
	x\plusone y:=x+y-(f(\Delta_{d-1}),2)
	\qquad\textrm{and}\qquad
	x\plustwo y:=x+y-(e_0+e_{d-1});
\]
These do not have neutral elements, as $(f(\Delta_{d-1}),2)$ and $e_0+e_{d-1}=(1,0,\dots,0,1)$
are not $f$-vectors of polytopes.
All three operations are equivalent to the usual vector addition,
via a suitable translation of~$\fset$ that sends
$f(\Delta_d)$, $(f(\Delta_{d-1}),2)$ resp.\ $e_0+e_{d-1}$ to the origin. 
If $x$ and $y$ satisfy the Euler equation,
and possibly also the Dehn--Sommerville equations, then
the same is also true for $x\pluszero y$, $x\plusone y$, and $x\plustwo y$.)  

We get the structure of an affine monoid $(\fset,\pluszero)$
\begin{compactitem}
	\item on the $f$-vector set $\fset=\fset(\mathcal{P}^2)$ of all $2$-polytopes,
	\item on the $f$-vector set of the $3$-polytopes $\fset(\mathcal{P}^3)$ 
		(by a 1906 lemma by Steinitz), as well as
	\item on the $f$-vector set $\fset(\mathcal{P}^d_s)$ of all simplicial $d$-polytopes,
	      for any $d\ge2$ (using the $g$-Theorem and addition of $M$-sequences).
	\item This \emph{fails} for general $4$-polytopes,
		as there are two $f$-vectors whose modified sum is not an $f$-vector,
	\item but we will prove an approximate version for the $f$-vector set of all $4$-polytopes.
\end{compactitem}

One would like to interpret the ``$\pluszero$'' addition operation
on $f$-vectors as the result of some glueing operation on polytopes that have the
$f$-vectors in question. However, we can describe such a glueing operation only on a restricted
class of ``adapter polytopes.''

For more general $d$-polytopes we discuss two natural glueing operations: The first one
is the glueing of two polytopes -- using a projective transformation if needed -- in a simplex facet.
However, 
\begin{compactitem}
	\item this operation is also not defined for all polytopes; we can define it
	on the class of $d$-polytopes with at least two simplex facets, such that the
	result of the glueing again lies in the same class;
	\item the corresponding $f$-vector operation is not given by $\pluszero$, 
		but rather by $x\plusone y:=x+y-(\Delta_{d-1},2)$;
	\item this addition operation cannot yield a monoid structure, as
    there is no neutral element. 
\end{compactitem}
With this, we get the structure of an affine semigroup $(\fset,\plusone)$
\begin{compactitem}
	\item on the $f$-vector set $\fset(\mathcal{P}^2)$ of all $2$-polytopes,
	\item on the $f$-vector set $\fset(\mathcal{P}^3)$ of all $3$-polytopes 
		 (since every $f$-vector is represented by a polytope with two triangle facets), and
	\item on the $f$-vector set $\fset(\mathcal{P}^d_s)$ of all simplicial $d$-polytopes
	     for all $d\ge2$ (where the semigroup addition is represented by glueing of polytopes).
	\item The semigroup addition $\plusone$ \emph{fails} for general $4$-polytopes,
		as there are two $f$-vectors whose modified sum is not an $f$-vector,
	\item but it is again available in an approximate version for general $4$-polytopes.
\end{compactitem}

The second geometric glueing operation is the 
``connected sum'' $P\#Q$ of two polytopes, which was apparently first used by Eckhoff in 1985.
However, 
\begin{compactitem}
	\item this operation is also not defined for all polytopes; we need extra conditions,
		e.g.\ that $P$ has a simple vertex and $Q$ has a simplex facet;
	\item the corresponding $f$-vector operation given by $x\plustwo y:=x+y-(e_0+e_{d-1})$;
	\item this addition operation also has no neutral element on $f$-vector sets of polytopes. 
\end{compactitem}
Nevertheless, we get the structure of an affine semigroup $(\fset,\plustwo)$
\begin{compactitem}
	\item on the $f$-vector set of the $2$-polytopes $\fset(\mathcal{P}^2)$ and
	\item on the $f$-vector set of the $3$-polytopes $\fset(\mathcal{P}^3)$.
	\item We also have an approximate semigroup structure on the
		 $f$-vectors of the $4$-polytopes, but
	\item Conjecture~\ref{conjecture:semigroup_d=4} says that 
			 the $f$-vector set of the $4$-polytopes $\fset(\mathcal{P}^4)$
			is closed under the operation $\plustwo$, so~that
			$(\fset(\mathcal{P}^4),\plustwo)$ is an affine semigroup.
\end{compactitem}

In this paper we proceed as follows:
In Section~\ref{sec:dim2and3} we describe the
		situation in low dimensions, $d\le3$.
Then in Section~\ref{sec:lattices} we describe the
affine lattices spanned by all $d$-polytopes,
and of all simplicial (equivalently: simple) $d$-polytopes, for all $d\ge2$:
As this does not seem to appear anywhere in the literature, we present it in detail
and with elementary, complete, uniform proofs.

The $f$-vector sets of the $d$-polytopes are contained 
in cancellative affine sub-monoids of the $f$-vector lattices that
also can be explicitly described,
see Section~\ref{sec:monoid_embeddings}.
The $f$-vectors of simplicial/\allowbreak simple
$d$-polytopes form a monoid by themselves.
This result does not extend canonically to the
$f$-vectors of general $d$-polytopes, as shown in Example~\ref{example:non-additive}.

In Section~\ref{sec:semigroups_by_glueing} we discuss situations
in which a monoid or semigroup structure can be represented geometrically,
by one of three types of ``glueing of polytopes.''  

Finally, in Section~\ref{sec:approximate_monoids} we show that
the $f$-vector set of the $4$-polytopes is an approximate semigroup:
Theorem \ref{thm:approximate_semigroup_d=4} says that
the sum of any two $f$-vectors of $4$-polytopes is close to an $f$-vector,
which can be realized by ``glueing the polytopes after modification.''

\section{Dimensions 2 and 3}%
\label{sec:dim2and3}%

For dimension $2$ we clearly have
\[
   \fset(\mathcal{P}^2) = \{ (n,n) : n\ge3 \}.
\]
This may be viewed as an (additive, commutative, cancellative) semigroup
with the usual vector addition:
Any sum of two $f$-vectors is an $f$-vector.
Moreover, it may be viewed as 
the set of all the integer points in a $1$-dimensional cone,
whose apex lies at $f(\Delta_3)=(3,3)$.
Thus with the addition
$(n,n)\pluszero(m,m)=(n+m-3,n+m-3)$ we have an affine monoid.
\smallskip

A simple and complete description of the
set of $f$-vectors of the $3$-dimensional polytopes was given by Steinitz in 1906:

\begin{lemma}[Steinitz's Lemma \cite{Stei3}]\label{lemma:Steinitz}
	The $f$-vectors of $3$-polytopes are the integer points in a $2$-dimensional affine cone
	whose apex lies at $f(\Delta_3)=(4,6,4)$\textup:
\[
   \fset(\mathcal{P}^3) = \{ (f_0,f_1,f_2)\in\Z^3 : 
   f_1=f_0+f_2-2,\ f_2\le2f_0-4,\ f_0\le2f_2-4\}.
\]
\end{lemma}

Let us note that this cone lies in the 
\emph{Euler hyperplane} $\EE_3=\{(f_0,f_1,f_2)\in\R^3:f_0-f_1+f_2=2\}$,
which is an affine plane in $\R^3$ that does not contain the origin.
Thus, in particular, a sum of two $f$-vectors is never an $f$-vector,
but it is close. However, translation of the $f$-vector set 
by the vector $-f(\Delta_3)$, or equivalently use of the addition
\[
   (x_0,x_1,x_2)\pluszero(y_0,y_1,y_2)=(x_0+y_0-4,x_1+y_1-6,x_2+y_2-4)
\]
yields the perfect monoid structure.

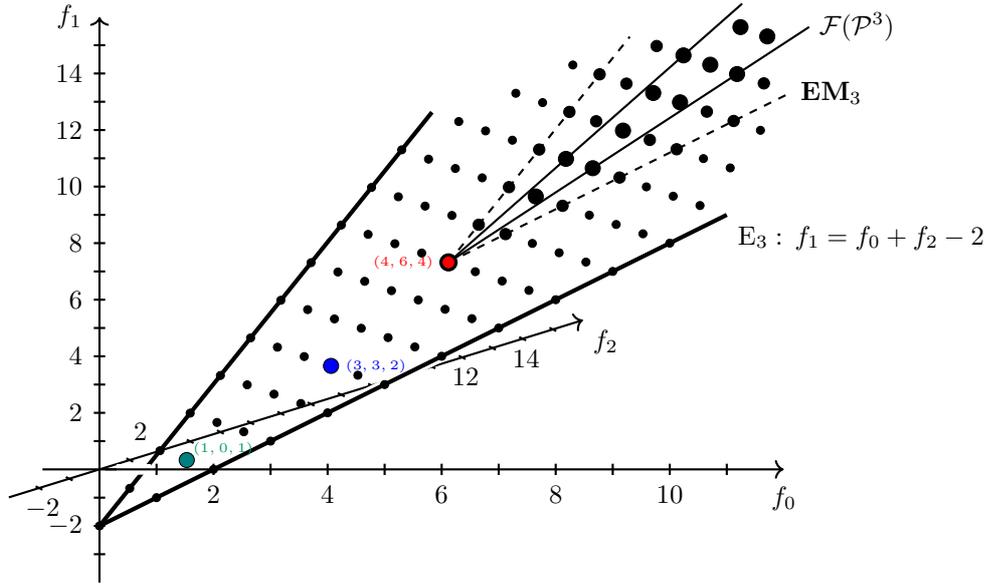
\begin{figure}[ht!]
	\begin{center}


\tikzstyle{Punkt} =          [draw, shape=circle, fill,inner sep=0pt,minimum size=3pt]
\tikzstyle{MittlererPunkt} = [draw, shape=circle, color=black,inner sep=1.5pt, minimum size=2pt, fill]
\tikzstyle{GrosserPunkt} =   [draw, shape=circle, color=black,inner sep=2pt, minimum size=3pt, fill]
\tikzstyle{RoterPunkt} = [draw,very thick, shape=circle,inner sep=2pt, minimum size=3pt, color=black, fill=red]
\tikzstyle{BlauerPunkt} = [draw, shape=circle,inner sep=2pt, minimum size=2pt, color=black, fill=blue]
\tikzstyle{GruenerPunkt} = [draw, shape=circle,inner sep=2pt, minimum size=2pt, color=black, fill=green!50!blue]

%

\begin{tikzpicture}[scale=.75,yscale=0.5, z={(0.53,0.33)},line join=round]



\draw[thick, ->] (-1,0,0) -- (12,0,0) node[below=0.1] {$f_0$};
\draw[thick, ->] (0,-4,0) -- (0,16,0) node[left=0.1] {$f_1$};
\draw[thick, ->] (0,0,-3) -- (0,0,16) node[below right] {$f_2$};

\foreach \x in {1,...,11}
\draw[thick] (\x,-0.2,0) -- (\x,0.2,0){};

\foreach \x in {2,4,...,10}
\node[scale=1.0, below=0.1] at (\x,0) {$\x$};

\foreach \y in {-3,-2,-1,1,2,...,15}
\draw[thick] (-0.1,\y,0) -- (0.1,\y,0) {};

\foreach \y in {-2,2,4,...,14}
\node[scale=1.0, left=0.1] at (0,\y) {$\y$};

\foreach \z in {-2,-1,1,2,...,15}
\draw[thick] (-0.06,0.06,\z) -- (0.06,-0.06,\z) {};

\foreach \z in {2}
\node[above left] at (-0.06,0.06,\z) {$\z$};
\foreach \z in {-2,12,14} 
\node[below] at (0.06,-0.06,\z) {$\z$};


\draw[white, line width=4pt] (4,2,0) -- (6,4,0);  
\draw[ultra thick] (0,-2,0) -- (11,9,0) node[below right]{$\EE_3:\,f_1=f_0+f_2-2$};
\draw[ultra thick] (0,-2,0) -- (0,9,11);
\draw[white, line width=4pt] (0.3,0,0) -- (0.9,0,0);
\draw[black, thick] (0.1,0,0) -- (0.9,0,0);

\foreach \i in {0,...,10}
   \foreach \j in {0,...,10}{
   	\pgfmathparse{max(\i+\j-2+\j,\j+\i-2+\i)}
      \ifdim \pgfmathresult pt < 22 pt \node[Punkt] at (\i,\i+\j-2,\j) {}; \fi
  }
\foreach \i in {4,...,10}
    \foreach \j in {4,...,10}{
    	\pgfmathparse{max(\i+\j-2+\j,\j+\i-2+\i)}
       \ifdim \pgfmathresult pt < 22 pt \node[MittlererPunkt] at (\i,\i+\j-2,\j) {}; \fi
   }
   

\draw[thick] (4,6,4) --  (6.5,13.5,9);
\draw[thick] (4,6,4) --  (9,13.5,6.5) node[right]{$\fset(\mathcal{P}^3)$};

\draw[thick,dashed] (4,6,4) --  (4,12,10);
\draw[thick,dashed] (4,6,4) -- (10,12,4) node[right]{$\EEM_3$};

\foreach \i in {4,...,8}
\foreach \j in {4,...,8}{
	\pgfmathparse{max(\j-\i,\i-\j)}
	\edef\temp{\pgfmathresult}
	\pgfmathparse{min(\i-4+1,\j-4+1)}
	\ifdim \temp pt < \pgfmathresult pt 
	  \pgfmathparse{\i+\j}
	  \ifdim \pgfmathresult pt < 16 pt \node[GrosserPunkt] at (\i,\i+\j-2,\j) {}; \fi
	\fi
}


\node[RoterPunkt] at (4,6,4) {};
\node at (4,6,4) [left=2pt] {\color{red}\tiny$(4,6,4)$};
\node[BlauerPunkt] at (3,3,2) {};
\node at (3,3,2) [right=2pt] {\color{blue}\tiny$(3,3,2)$};
\node[GruenerPunkt] at (1,0,1) {};
\node at (.9,0.4,1) [right=1pt] {\color{ForestGreen}\tiny$(1,0,1)$};

\end{tikzpicture}

	\end{center}
	\caption{\label{fig:planeE3}%
		Integer points on the Euler hyperplane $\EE_3$ with $f_0,f_2\ge0$.
		The subset of lattice points with $f_0,f_2\ge4$ forms the Euler monoid $\EEM_3$,
		whose neutral element is $(4,6,4)$.
		This in turn contains as subset the $f$-vectors of $3$-polytopes,
		as described by Steinitz's Lemma~\ref{lemma:Steinitz},
		which again is given by all lattice points in a cone with apex $(4,6,4)$.
		Also shown are two other ``base points'' $(3,3,2)$ and $(1,0,1)$
		on the Euler hyperplane, will be discussed later.}
\end{figure}

The cone described by Steinitz' lemma is spanned by the $f$-vectors 
of the bipyramid $(5,9,6)$ 
and of the prism $(6,9,5)$. 
The Hilbert basis for the cone, that is, the minimal set of generators
for the monoid, additionally contains the $f$-vector 
of the square pyramid $(5,8,5)$.

Before we proceed with monoid structures, we will here
discuss the lattices spanned by $f$-vector sets (where you might
take the cases $d=2,3$ as examples.)

\section{Lattices spanned by \emph{f}-vector sets}%
\label{sec:lattices}%

A \emph{lattice} $\Lambda\subset\R^d$ is a discrete abelian subgroup,
and an \emph{affine lattice} is any translate of a lattice.
We will exclusively deal with \emph{integer lattices}, which are
subsets of~$\Z^d$.
The \emph{affine lattice} $\aff_\Z S$ spanned by a subset $S\subseteq\Z^d$
is given by all the integral affine combinations
$a_1v_1+\dots+a_nv_n$ with $n\ge1$, $v_1,\dots,v_n\in S$, and
$a_1,\dots,a_n\in\Z$ with $a_1+\dots+a_n=1$.

In this section, we prove that
\begin{compactitem}
	\item the affine lattice $\aff_\Z\fset(\mathcal{P}^d)$ spanned by the
		$f$-vectors of $d$-polytopes contains \emph{all}
		integer points on the Euler hyperplane $\EE_d$, and
	\item the affine lattice $\aff_\Z\fset(\mathcal{P}^d_s)$ spanned by the
		$f$-vectors of simplicial $d$-polytopes contains \emph{all}
		integer points on the Dehn--Sommerville subspace $\DS_d$.
\end{compactitem}
The first result we obtain by extending the proof given by
Walter Höhn \cite[pp.~29,30]{Hoehn-DissETH1953}
when in his 1953 PhD thesis at ETH Zurich he first established
that $\aff_\R\fset(\mathcal{P}^d)=\EE_d$.
 
The second one is based on a proof given by
Victor Klee in 1964 \cite[Prop.~3.4]{Klee3}, who used joins of 
simplex boundaries (that is, the boundary of a direct sum of simplices) 
in order to establish that 
$\dim\lin_\R\fset(\mathcal{P}^d_s)\ge\lfloor\tfrac12(d+1)\rfloor$.
Alternatively, the result for simplicial polytopes can be derived from Björner's 
\cite{Bjo2} \cite{Bjo6} matrix formulation of the $g$-Theorem --
without using any non-trivial part of the proof of the $g$-Theorem.
(This method will be employed later in this paper.)

We here give a self-contained treatment, which also emphasizes
that both results can be obtained with the same type of argument, applied to
different sets of $d$-polytopes with at most $d+2$ vertices.	

Curiously enough, the question about the affine lattice of $f$-vectors,
i.e.\ whether there are hidden parity or congruence constraints,
for general or for simple/simplicial polytopes,
apparently had not been asked or answered up to now for general or for
simplicial $d$-polytopes. 
However, the question was explicitly treated for cubical polytopes,
where Babson \& Chan \cite{BabsonChan3} obtained 
partial answers motivated by the result by Blind \& Blind \cite{BlBl1} 
that all even-dimensional cubical polytopes have an even number of vertices.
	
\subsection{The \emph{f}-vectors of \emph{d}-polytopes with at most \emph{d}+2 vertices}%
\label{subsec:d+2}%

It is not hard to see \cite[Sect.~6.1]{Gr1-2} 
that any combinatorial type of a $d$-polytope
with at most $d+2$ vertices
arises as an $m$-fold pyramid over a 
direct sum of a $k$-simplex and an $\ell$-simplex,
\[
	P_{k,\ell}^m := 
	\Delta_{m-1}*(\Delta_k\oplus\Delta_\ell)
	\qquad\textrm{for } k,\ell,m\ge0.  
\]
These are polytopes of dimension $d:=k+\ell+m$ with $f_0=k+\ell+m+2=d+2$ vertices, 
except that for $k=0$ we get $P_{0,\ell}^m=\Delta_{m-1}*\Delta_\ell=\Delta_{\ell+m}$,
which has $\ell+m+1=d+1$ vertices.
(To obtain \emph{distinct} combinatorial types, we
could put the additional conditions that $k\le\ell$ 
and that $m=0$ if $k=0$. Indeed, $P_{0,\ell}^m$
is an $(\ell+m)$-simplex.)

\begin{lemma}\label{lemma:f-klm}
	For $k,\ell,m\ge0$, the $f$-vector of $P_{k,\ell}^m$ is given by
	\[
	f_i(P_{k,\ell}^m) =
	  \binom{k+\ell+m+2}{i+1}
	- \binom{\ell+m+1}{i-k}
	- \binom{k+m+1}{i-\ell}
	+ \binom{m+1}{i-k-\ell}
	\qquad\textrm{for } -1\le i\le k+\ell+m.
	\]
\end{lemma}

\begin{proof}
	We describe the faces (and their vertex sets) explicitly,
	in terms of the vertex sets 
	$V_k:=V(\Delta_k)$, $V(\Delta_\ell)$, and $V(\Delta_{m-1})$
	of the three simplices that $P_{k,\ell}^m$ is constructed from.
	
	As the direct sum $\Delta_k\oplus\Delta_\ell$ is a simplicial polytope,
	the $i$-faces that do not contain $\Delta_k\oplus\Delta_\ell$
	are all simplices, so they have $i+1$ vertices. Their vertex sets are given by
	all the $\binom{k+\ell+m+2}{i+1}$  $(i+1)$-subsets of vertices of $P_{k,\ell}^m$,
	minus those $\binom{\ell+m+1}{(i+1)-(k+1)}$ subsets that contain $V(\Delta_k)$ and
	minus those $\binom{k+m+1}{(i+1)-(\ell+1)}$ subsets that contain $V(\Delta_\ell)$, 
	plus  those $\binom{m}{(i+1)-(k+\ell+2)}$ subsets that we have just subtracted twice,
	as they contain $V(\Delta_k)$ \emph{and} $V(\Delta_\ell)$.

	All the other $i$-faces are multiple pyramids over the $(k+\ell)$-face
    $\Delta_k\oplus\Delta_\ell$, and these are counted by~$\binom{m}{i-(k+\ell)}$.
	Thus we obtain
	\[
	f_i(P_{k,\ell}^m) =
	  \binom{k+\ell+m+2}{i+1}
	- \binom{\ell+m+1}{i-k}
	- \binom{k+m+1}{i-\ell}
	+ \binom{m}{i-k-\ell-1}
	+ \binom{m}{i-k-\ell}.
	\] 
	Combining the last two terms completes the proof.
\end{proof}

\subsection{The lattice of \emph{f}-vectors of all \emph{d}-polytopes}

We first treat the $f$-vector sets of general $d$-polytopes.
These lie on the \emph{Euler hyperplane}
\[
	\EE_d := \{(f_0,\dots,f_{d-1}): f_0-f_1+\dots+(-1)^{d-1}f_{d-1}=1-(-1)^d\} \subset\R^d.
\]
We will prove that \emph{all} integer points on this
hyperplane are affine combinations of $f$-vectors of $d$-polytopes.

\begin{theorem}
	For any $d\ge2$, we have 
	\[
	  \aff_\Z(\fset(\mathcal{P}^d)) = \EE_d\cap\Z^d.
	\]
\end{theorem}

Thus all the parity constraints, and other congruences that are valid for
the $f$-vectors of $d$-polytopes, are implied by Euler's equation.
For example, the equation $f_0-f_1+f_2=2$ for $3$-polytopes implies that
the number of faces $f_0+f_1+f_2$ is even for every $3$-polytope.

\begin{proof}
It will be convenient to use \emph{extended $f$-vectors}
\[
	\hat{f}(P):=(1,f(P))=(f_{-1},f_0,f_1,\dots,f_{d-1})\in\Z^{d+1}.
\]
which are collected in the extended $f$-vector set 
$\hat\fset(\mathcal{P})\subset\Z^{d+1}$.

To prove the theorem we proceed as follows:
We exhibit $d$ $d$-polytopes, 
form a matrix $N_d\in\Z^{d\times(d+1)}$ that 
has their extended $f$-vectors as the rows,
and show that this matrix has full rank $d$, and that all integer vectors in the row space of~$N_d$
are integer combinations of the rows of the matrix. 

Let's do it! We follow Höhn \cite[pp.~29,30]{Hoehn-DissETH1953}
to use the polytopes 
	\[ 
	  P_k(d):= P^{d-k-1}_{k,1} =
	  \Delta_{d-k-2}*(\Delta_k\oplus\Delta_1)
	  \qquad\textrm{for }k=0,1,\dots,d-1.
	\]
	Höhn describes them as being glued together from two $d$-simplices;
	we have presented them as multiple pyramids over the bipyramid over a simplex
	$\pyr^{d-k-1}\bip(\Delta_k)$.
Lemma~\ref{lemma:f-klm} yields the extended $f$-vectors:
	\begin{equation}\label{eq:fi-Pkd}
	f_i(P_k(d)) = \binom{d+1}{i+1} + \binom{d}{i} - \binom{d-k}{d-i}.
	\end{equation} 
Now we build the matrix 
\[
\widetilde{N}_d := 
\Big(
  f_i\big(P_k(d)\big) 
\Big)_{\substack{0\le k\le d-1\\
   \hspace{-5pt}-1\le i\le d-1}}
  \in \Z^{d\times(d+1)},
\]
which has the extended $f$-vectors of the $d$ polytopes as rows.
As the first two terms in (\ref{eq:fi-Pkd}) do not depend on $k$,
we see that after  
subtracting from each row the row above it (proceeding from the bottom
row towards the top) we arrive at an upper-triangular matrix ${N}_d$,
with entries
\[
- \binom{d-k}{d-i} + \binom{d-k+1}{d-i} = \binom{d-k}{d-i-1}  
\]
in all rows below the top row, 
which \emph{obviously} integrally spans the full lattice in its row space.
(The rows below the top row, whose first nonzero entry is~$+1$,
also form part of the Pascal triangle.)
For example, for $d=3$ we obtain
\[
  \widetilde{N}_3 = 	\left(
  		\begin{array}{llll}
		1 & 4 & 6 & 4\\
		1 & 5 & 8 & 5\\
		1 & 5 & 9 & 6\\
		\end{array} 
		\right)
		\quad\rightsquigarrow\quad
  {N}_3 =
	 	\left(
		\begin{array}{llll}
		1 & 4 & 6 & 4\\
		0 & 1 & 2 & 1\\
		0 & 0 & 1 & 1\\
		\end{array}
		\right).
\]
\vskip-6mm%
\end{proof}

\subsection{The lattice of \emph{f}-vectors of simplicial \emph{d}-polytopes}

Now we treat the $f$-vectors of the (combinatorial types of)
\emph{simplicial} $d$-polytopes. (By duality, this yields the corresponding
results for the simple $d$-polytopes.)
For simplicial polytopes the Dehn--Sommerville equations
-- additional linear equations satisfied by the $f$-vectors of simplicial polytopes,
see \cite[Sect.~9.2]{Gr1-2} --
imply that the $f$-vectors lie in an affine subspace of~$\R^d$
of dimension $\dim \DS_d=\dhalfdown$, 
the \emph{Dehn--Sommerville subspace}
\[
	\DS_d := \aff_\R(\fset(\mathcal{P}^d_s)).
\]

Again, we claim that \emph{all} the integer points on that subspace
are integral affine combinations of $f$-vectors of simplicial $d$-polytopes.

\begin{theorem}\label{thm:DS-lattice}
	For any $d\ge2$, we have 
	\[
	  \aff_\Z(\fset(\mathcal{P}^d_s)) = \DS_d\cap\Z^d.
	\]
\end{theorem}

Thus all the parity constraints, and all other congruences
that are valid for
the $f$-vectors of simplicial polytopes, are implied by the Dehn--Sommerville equations.
For example, the equation $3f_2=2f_1$ for simplicial $3$-polytopes implies that
$f_2$ is even and $f_1$ is divisible by $3$. Less trivial congruences are analyzed
in Björner \& Linusson \cite{BjornerLinusson}.

\begin{proof}
It will again be convenient to use {extended $f$-vectors}.
To prove the theorem we now exhibit $\lfloor d/2\rfloor+1$ simplicial $d$-polytopes, 
form a matrix 
$M_d\in\Z^{(\dhalfdown+1)\times(d+1)}$ 
that has their extended $f$-vectors as the rows,
and show that $M_d$ has full row rank $\dhalfdown+1$ and that all integer vectors in the row space of~$M_d$
are integer combinations of the rows of the matrix. 

Let's do it!  
For the polytopes we take the polytopes 
	\[
	  S_k(d):= P^{0}_{k,d-k} =%
	  \Delta_k\oplus\Delta_{d-k}
	  \qquad\textrm{for }k=0,1,\dots,\dhalfdown.
	\]
Lemma~\ref{lemma:f-klm} yields the extended $f$-vectors:
	\begin{equation}\label{eq:fi-Skd}
	f_i(S_k(d)) = \binom{d+2}{i+1} - \binom{k+1}{i-d+k} - \binom{d-k+1}{i-k}
	\qquad \textrm{for }-1\le i\le d-1.
	\end{equation} 
(Note that if we evaluate this for $k=-1$ and $-1\le i\le d-1$, we get a zero vector.)
Now we build the matrix 
\[
\widetilde{M}_d := 
\Big(
  f_i\big(S_k(d)\big)
\Big)_{\substack{0\le k\le \dhalfdown\\
   \hspace{-5pt}-1\le i\le d-1}}
  \in \Z^{(\dhalfdown+1)\times(d+1)},
\]
which has the extended $f$-vectors of the $\dhalfdown+1$ polytopes as rows.
The first two terms in (\ref{eq:fi-Skd}) do not depend on $k$.
After subtracting from each row the row above it (proceeding from the bottom
row towards the top) we arrive at an ``upper-triangular'' $(\dhalfdown+1)\times(d+1)$-matrix $M_d$
with $(k,i)$-entries
\[
- \binom{k}{i-d+k} + \binom{d-k+1}{i-k+1},
\]
which \emph{obviously} integrally spans the full lattice in its row space.
We note that the matrix $M_d$ is exactly the ``McMullen matrix''
from Björner's \cite{Bjo2} \cite{Bjo6} matrix formulation of the $g$-Theorem, 
see \cite[Thm.~8.36]{Z35}.
For example, for $d=4$ we obtain
\[
  \widetilde{M}_4 = \left(
  		\begin{array}{rrrrr}
		1 & 5 & 10 & 10 & 5\\
		1 & 6 & 14 & 16 & 8\\
		1 & 6 & 15 & 18 & 9\\
		\end{array} 
		\right)
		\quad\rightsquigarrow\quad
	{M}_4 =
	 	\left(
		\begin{array}{rrrrr}
		1 & 5 & 10 & 10 & 5\\
		0 & 1 &  4 &  6 & 3\\
		0 & 0 &  1 &  2 & 1\\
		\end{array}
		\right).
\]
\vskip-6mm%
\end{proof}

\section{Cancellative affine monoids containing \emph{f}-vector sets}%
\label{sec:monoid_embeddings}%

The lattices $\EE_d$ and $\DS_d$ give ``outer descriptions''
of the $f$-vector sets of all $d$-polytopes, 
and of all simplicial $d$-polytopes, respectively.
We now head for tighter descriptions.

For this we again start with the case of general $d$-polytopes: 
The set $\fset(\mathcal{P}^d)$ is contained
in a monoid of all integer points in a pointed $(d-1)$-dimensional 
cone, which we call the ``Euler monoid.''
We describe the combinatorics of the cone, and also the Hilbert basis.

After this, we describe the ``Dehn--Sommerville monoid'' that
contains the $f$-vectors of simplicial $d$-polytopes -- where we
get much more, namely an exact monoid structure on the set of $f$-vectors.

\subsection{The Euler monoid}

The \emph{trivial lower bound inequality} $f_i(P)\ge f_i(\Delta_d)=\binom{d+1}{i+1}$,
or in vector notation with componentwise comparison 
$f(P)\ge f(\Delta_d)$, is valid for all polytopes;
see Grünbaum \cite[Ex.~3.1.8, p.~36]{Gr1-2}.

\begin{definition}[The Euler monoid]
	The \emph{Euler monoid} $\EEM_d=(\mathrm{EM}_d,\pluszero)$ is given by the 
	set of integer points $x\in\Z^d$ that satisfy the Euler equation
	and the trivial lower bound inequality,
	\[
	\mathrm{EM}_d := \{ x\in\EE_d: x\ge f(\Delta_d)\} \cap \Z^d,
	\]
	with the reduced addition operation
	\[
	x \pluszero y := x+y-f(\Delta_d).
	\] 
\end{definition}

Due to the trivial lower bound inequality,
the $f$-vector set $\fset(\mathcal{P}^d)$ is contained
in the ground set of the Euler monoid,
\[
	\fset(\mathcal{P}^d) \subset \mathrm{EM}_d.
\]
This embeds the $f$-vector set into an algebraic structure that we fully understand.
(Cf.~Figure~\ref{fig:planeE3}.)

\begin{proposition}[The Euler monoid]
	The Euler monoid consists of all integer points in the
	pointed polyhedral cone 
	\[
	\mathrm{CEM}_d := \{ x\in\EE_d: x\ge f(\Delta_d)\},
	\]
	which has the combinatorics of a cone over a product of two simplices
	$\Delta_{\dhalfup-1}\times\Delta_{\dhalfdown-1}$.
	
	The minimal set of generators for the monoid,
	that is, the reduced Hilbert basis for the cone, is
	given by the vectors $f(\Delta_d)+e_i+e_j$
	with $i,j\in\{0,1,\dots,d\}$, where $i$ is even and $j$ is odd.
\end{proposition}

\begin{proof}
	One easily verifies that the 
	$\dhalfup+\dhalfdown=d$ facets of this cone are given
	by the $d$ inequalities $x_i\ge f_i(\Delta_d)$.
	The $(\dhalfup+1)(\dhalfdown+1)$ extremal rays are given by the 
	direction vectors
	$e_i+e_j$ with $i,j\in\{0,1,\dots,d\}$, where $i$ is even and $j$ is odd.
	 
	Any vector $x$ in the cone has coordinates $x_i\ge f_i(\Delta_d)$.
	If $x$ is an integer point in the cone that is 
	not the apex, then due to the Euler equation it needs to have
	both an even index $i$ and an odd index $j$ such that
	$x_i$ and $x_j$ are larger than
	the lower bound. Due to integrality we get that $x-(e_i+e_j)$ lies in the cone.
	By induction, this verifies the Hilbert basis property.
\end{proof}

The inequalities 
$f_i(P)\ge f_i(\Delta_d)=\binom{d+1}{i+1}$
that we used in the construction of the Euler monoids are very weak.
Tighter descriptions are available, e.g.\ for any $d$
we could use that $f_1\ge\frac d2f_0$ and $f_{d-2}\ge\frac d2f_{d-1}$.

Ideally, one would want to use a tight description
of the closed convex cone with apex $f(\Delta_d)$
spanned by the $f$-vectors of $d$-polytopes.
However, such a description may not be finite, and it is 
not available even for the first non-trivial
case $d=4$: Here
five tight linear inequalities are known, the
two given by $f_0,f_3\ge0$, the two we just named 
($f_1\ge2f_0$ and $f_2\ge2f_3$), and a lower bound inequality 
that in Ziegler \cite{Z82} we interpreted as the
``fatness lower bound'' $F=\frac{f_1+f_2-20}{f_0+f_3-10}\ge\frac52$.
A sixth inequality, which would give an upper bound on $F$,
possibly $F\le9$, is missing.
Our knowledge is even less complete in higher dimensions;
see e.g.\ Björner \cite{Bjo2}, Eckhoff \cite{eckhoff:_combin2}, and Werner \cite[Chap.~6]{Werner-diss}.
\smallskip

At this point we leave the topic of ``tighter outer descriptions''
and instead look for the ``inner algebraic structure''
of the $f$-vector set $\fset(\mathcal{P}^d)$.

Indeed, we might be a bit bold and wonder whether possibly
$(\fset(\mathcal{P}^d),\pluszero)$ is a monoid by itself,
that is, whether for any $f$-vectors $x,y$ the reduced sum
$x\pluszero y$ is again an $f$-vector.
This is true for $d=2$ and $d=3$,
as we have seen in Section~\ref{sec:dim2and3}.

For $d=4$, this may still be true ``most of the time,''
say with only finitely many exceptions.
However, the following example -- the only one we know of,
except for a second one that arises from duality --
shows that if $P,Q$ are $4$-polytopes, 
then $f(P)+f(Q)-f(\Delta_d)$ is \emph{not in general} an $f$-vector
of a $4$-polytope.

\begin{example}\label{example:non-additive}
	Take $P=(P^0_{1,3})^*=\Delta_3\times\Delta_1=\prism(\Delta_3)$ and
	$Q=P^0_{2,2}=\Delta_2\oplus\Delta_2$,
	with $f$-vectors
	$f(P)=(8, 16, 14, 6)$ and
	$f(Q)=(6, 15, 18, 9)$.
	In this case there is no $4$-polytope with $f$-vector 
	$f(P)\pluszero f(Q)=(9, 21, 22, 10)$, 
	as this violates the fatness lower bound inequality for $4$-polytopes
	$F(P)=\frac{f_1+f_2-20}{f_0+f_3-10}\ge\frac52$, see 	
	Ziegler \cite[Sects.~4 and 5]{Z82}.
	
	(The fatness lower bound $F(P)\ge\frac52$ follows from the ``$g_2^{tor}(P)\ge0$'' inequality
	$f_{03} - 3(f_0 + f_3) \ge -10$, which Kalai \cite{kalai87:_rigid_i} derived from
	rigidity theory, combined with the
	``center boolean'' lower bound $f_{013}\ge 3f_{03}$,  
	using some generalized Dehn--Sommerville equations from Bayer \& Billera \cite{BaBi}.) 
	According to Brinkmann \& Ziegler \cite[Table~1]{Z157}, there is also no
	cellular sphere with the $f$-vector in question.)
\end{example}

How was this Example \ref{example:non-additive} found?
Combining the $f$-vectors of the $4$-polytopes with
at most $7$ vertices from Perles's analysis in \cite[Sect.~6.1]{Gr1-2} 
and with $8$ vertices by Altshuler \& Steinberg in \cite{altshuler84} and \cite{altshuler85} with the classification of $f$-vectors with
$9\le f_0,f_3$ and $f_0+f_3\le22$ 
in \cite[Table~1]{Z157} we obtained a complete list of
all $184$ $f$-vectors of $4$-polytopes with $f_0+f_3\le22$.
Moritz Firsching determined that, up to duality,
the above pair is the only one such that $f(P)\pluszero f(Q)$ satisfies $f_0+f_3\le22$,
but is not an $f$-vector.

Let us also note that $f(P)\pluszero f(Q)=(9, 21, 22, 10)$
is ``close to an $f$-vector'': For example, there
are $4$-polytopes with $f$-vector $(9, 22, 23, 10)$.
This observation will be followed-up in Section~\ref{sec:approximate_monoids}.

\subsection{The Dehn--Sommerville monoid}

For simplicial $d$-polytopes, however, we get a perfect 
semigroup structure.

\begin{definition}[The Dehn--Sommerville monoid]
	For any $d\ge2$,
	the \emph{Dehn--Sommerville monoid} is the
	pair $\DSM_d:=(\fset(\mathcal{P}^d_s),\pluszero)$.
\end{definition}

\begin{theorem}
	For any $d\ge2$, $\DSM_d$ is a cancellative affine monoid.
\end{theorem}

\begin{proof}
	For this we rely on the $g$-Theorem \cite[Thm.~8.36]{Z35}, 
	which in Björner's \cite{Bjo2} \cite{Bjo6} matrix formulation says that
	\begin{quote}
	\emph{the extended $f$-vectors $\hat{f}(P)$ of simplicial $d$-polytopes
		are exactly \\
		the vector-matrix products $(1,x)M_d$},
	\end{quote} 
	where
	\begin{compactitem}[ -- ]
		\item $M_d\in\Z^{(\dhalfdown+1)\times(d+1)}$ 
			is the matrix that we have already met
			in the proof of Theorem~\ref{thm:DS-lattice},
			whose first row is the extended $f$-vector $\hat{f}(\Delta_d)$, and
		\item  $(1,x)\in\Z^{1+\dhalfdown}$
			is an $M$-sequence, that is, by Macaulay's theorem
			(see e.g.\ \cite[Thm.~8.34]{Z35}),
			the $f$-vector of a multicomplex.
	\end{compactitem} 
	Now if $(1,x)$ and $(1,y)$ are $f$-vectors of multicomplexes,
	then so is $(1,x+y)$, as we may assume that the multicomplexes
	corresponding to $(1,x)$ and $(1,y)$ have disjoint vertex sets.
	
	Thus if $(1,x)M_d$ and $(1,y)M_d$ are $f$-vectors of simplicial 
	$d$-polytopes, then so is
	\[
		(1,x)M_d \pluszero (1,y)M_d 
		= (1,x + y)M_d,
	\] 
	and we are done.	
\end{proof}

\begin{remark}
	$\DSM_d:=(\fset(\mathcal{P}^d_s),\pluszero)$ is not 
	an \emph{affine semigroup} in the sense defined
	by Miller \& Sturmfels \cite[Thm.~7.4]{miller04:_combin_commut_algeb},
	as it fails an explicit assumption made there (and elsewhere):
	For $d\ge4$ this semigroup is \emph{not} finitely-generated.

	This is due to the fact that some constraints defining
	$M$-sequences, such as $g_2\le\binom{g_1+1}2$, 
	and the corresponding $f$-vector inequalities like $f_1\le\binom{f_0}2$,
	are \emph{concave} rather than convex.
	This leads us to 
	$M$-sequences like $(1,t,\tfrac12t(t+1),0\dots,0)$
	for $t\ge1$, which cannot be written as sums of other
	nonzero $M$-sequences -- and similarly for the corresponding
	$f$-vectors of polytopes.
\end{remark}

\section{Semigroup structures via glueing}%
\label{sec:semigroups_by_glueing}%

It is natural to ask, once we know that
$f(P)\pluszero f(Q)$ is the $f$-vector $f(R)$ of a polytope,
whether $R$ may be obtained by ``glueing'' $P$ and $Q$ in some way.

In general, there seems to be no easy positive answer to this,
even if we are permitted to replace $P$, $Q$, and $R$ by
other polytopes  $P'$, $Q'$, and $R'$ from the same class that have the same $f$-vector:
For all natural glueing operations to consider, 
there are restrictions on the polytopes that can be glued.

This allows for three different approaches, all of them of interest,
and all three of them will be briefly discussed in the following:
\begin{enumerate}
	\item Restrict to polytope subclasses which can be glued in order
		to represent the addition $f(P)\pluszero f(Q)$,
	\item modify the addition operation into a new one, like
		 $f(P)\plusone f(Q)$ or $f(P)\plustwo f(Q)$,
			whose output can be represented
			by a glueing operation like ``glueing in simplex facets'' 
			or the ``connected sum,'' 
			or
	\item modify the polytopes until they can be glued -- and thus get
		``approximate addition.''
\end{enumerate}

\subsection{Monoids of polytopes that can be glued}

In general, it is not clear what $d$-polytope would
have the $f$-vector $f(P)\pluszero f(Q)$. 

\begin{proposition}\label{prop:adapter-glueing}
	Let $P$and $Q$ be $d$-polytopes such that
	$P$ has a simple vertex $v$ that is contained in $d$ simplex facets, 
	and $Q$ also has a simplex facet $F$.
	Then one can delete the vertex $v$, resulting in a polytope $P'$
	with one vertex less and with a new simplex facet $F'$ spanned by the neighbors 
	of~$v$ in~$P$, and by a projective transformation $\pi$ (if necessary)
	one can position $\pi(P')$ and $Q$ such that 
	$\pi(P')\cap Q=\pi(F')=F$, 
	$\pi(P')\cup Q$ is convex,
	and all proper faces of $\pi(F')=F$ are also faces of $\pi(P')\cup Q$.

	In this situation, 
\[
	f(\pi(P')\cup Q) = f(P) \pluszero f(Q).
\]
\end{proposition}

\begin{figure}[ht!]
	\begin{center}
			\input{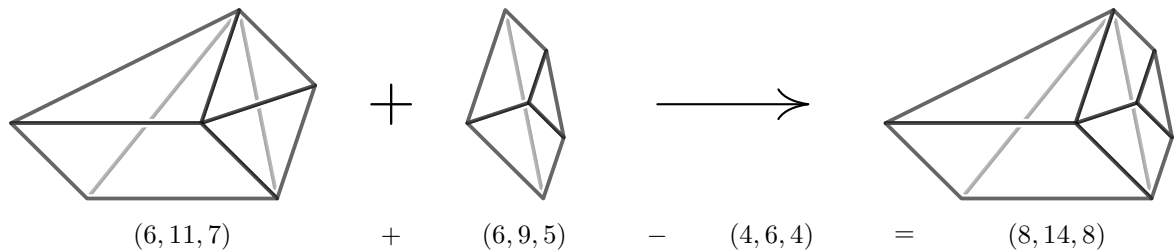} 
	\end{center}
	\caption{\label{fig:Polytope-a}%
		An example for the adapter glueing construction, and its effect on $f$-vectors, 
		for $d=3$.}
\end{figure}

Figure~\ref{fig:Polytope-a} illustrates this.

Note that in this situation (which is reminiscent of valuations!), 
the resulting polytope $\pi(P')\cup Q$ need not have
any of the properties that we required about $P$ or $Q$.
It may not have any simplex facet at all.
Thus, if we want to put the same assumptions on $P$, $Q$, and $\pi(P')\cup Q$,
we have to demand more.
The following class of $d$-polytopes ``does the job.''

\begin{definition}[The class of adapter polytopes]
	Let $\mathcal{P}^d_{ap}$, the \emph{adapter polytopes}, be the class of all $d$-polytopes
	$P$ that have 
		a simple vertex $v$ that is contained in $d$ simplex facets of $P$, and
		an additional simplex facet $F$ that does not contain $v$.
\end{definition}

\begin{proposition}
	The $f$-vectors of the $d$-dimensional adapter polytopes, with modified addition,
	form an affine monoid
	\[
	\big(\fset(\mathcal{P}^d_{ap}), \pluszero ).
	\]
	Moreover, the addition in this monoid is realized by
	\[
	f(P)\pluszero f(Q) = f(R),
	\]
	where $R=P'\cup Q$ is obtained by the glueing of Proposition~\ref{prop:adapter-glueing}.
\end{proposition}

One can and must ask, of course, whether the adapter polytopes
are very special, or quite general. Is~$\fset(\mathcal{P}^d_{ap})$
a good approximation to $\fset(\mathcal{P}^d)$?
Clearly we have $\mathcal{P}^d_{ap}\subset \mathcal{P}^d$:
The sets do not coincide, even for $d=3$, where 
the only polytope $P$ with $f$-vector $f(P)=(6,9,5)$ is a triangular prism,
which is not an adapter polytope.
Its dual, however, the triangular bipyramid, is an adapter polytope --
so the class of adapter polytopes is \emph{not} closed under duality.
Also the glueing operation is not specified completely,
and the combinatorial type of the result may depend on the order
of the polytopes, but the $f$-vector of the result is unique.

In Section~\ref{sec:approximate_monoids} we will show, however,
that for $d\le4$ every $f$-vector of a $d$-polytope is ``close''
to the $f$-vector of an adapter polytope.
(For every simplicial $d$-polytope, this is easy to see, for all $d\ge2$.)

\subsection{Semigroup structure induced by connected sums}

For more general $d$-polytopes we will consider two natural ways to combine two of them.

\begin{figure}[ht!]
	\begin{center}
			\input{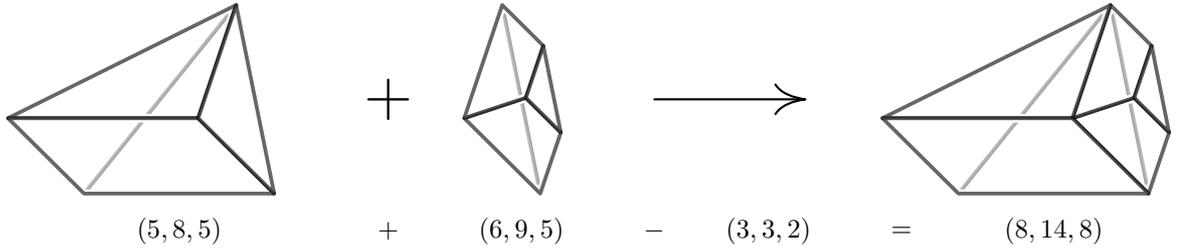} 
	\end{center}
	\caption{\label{fig:Polytope-b}%
	Examples for the ``glueing in simplex facets,''
	and its effect on $f$-vectors, for~$d=3$.}
\end{figure}

The first one is ``glueing in simplex facets,'' illustrated in the
Figure~\ref{fig:Polytope-b}:
For this let $P$ and $Q$ be $d$-polytopes which each have a simplex facet, $F\subset P$ and $F'\subset Q$.
After a suitable projective transformation $\pi$ we get
that $R:=\pi(P)\cup Q$ is a convex polytope such that all
proper faces of $\pi(P)\cap Q=\pi(F)=F'$ are also faces of $R$, and such that all proper faces
of $\pi(F)=F'$ are also faces of~$R$. In this case we have
\[
	f(R) = f(\pi(P)\cup Q) = f(P) + f(Q) - (f(\Delta_{d-1}),2) =: f(P) \plusone f(Q).
\] 
(Here $f(\Delta_{d-1})$ is a vector in $\Z^{d-1}$, so appending a ``$2$'' it we 
get a vector in $\Z^d$ that may be interpreted as the $f$-vector of a cellular $(d-1)$-sphere
that arises from two copies of $\Delta_{d-1}$ by identifying their boundaries with each other.)

\begin{figure}[ht!]
	\begin{center} 
			\input{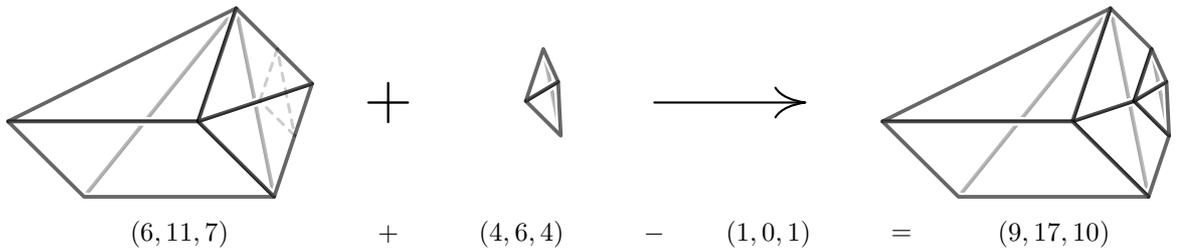}
	\end{center}
	\caption{\label{fig:Polytope-c}%
	An example for the the connected sum operation
	and its effect on $f$-vectors, for~$d=3$.}
\end{figure}

The second glueing operation, illustrated in Figure~\ref{fig:Polytope-c},
we consider is the 
``connected sum'' operation introduced by Eckhoff in 1985 \cite{eckhoff:_combin} 
\cite{eckhoff:_combin2}, compare \cite[p.~279]{Z35}:
Let $P$ be a $d$-polytope with a simple vertex and $Q$
a $d$-polytope with a simplex facet.
Then $P\#Q$ is formed by cutting off the simple vertex,
which generates a new simplex facet with $d$ new vertices, and then
glueing $Q$ into the newly formed facet, using a projective transformation if necessary. 
The process is quite similar to the one just described for adapter polytopes.
See Ziegler \cite[Example 8.41]{Z35} for details.
There is, however, one essential difference: 
For Eckhoff's glueing we get a different modified addition:
\[
	f(P\#Q) = f(P) + f(Q) - (1,0,\dots,0,1) =: f(P) \plustwo f(Q).
\]

Neither of these modified addition operations yields a monoid structure,
as there are no polytope with $f$-vectors $(f(\Delta_{d-1}),2)$ or $(1,0,\dots,0,1)$.
On the other hand, both of them yield
``geometrically realized'' semi-group structures on 
the $f$-vector sets of various classes of polytopes.

We record several instances of this. Their proofs are quite straightforward
and thus mostly omitted.

\begin{proposition}
	Let $\mathcal{P}^d_{2s}$ be the class of all $d$-dimensional
	polytopes with at least two simplex facets.
	Then the $f$-vector set forms an affine semigroup 
	\[
	\big(\fset(\mathcal{P}^d_{2s}),\plusone\big),
	\]
	where the modified addition can be represented by 
	``glueing in simplex facets'': 
	The result $\pi(P)\cup Q$ of glueing any two polytopes $P,Q\in\mathcal{P}^d_{2s}$
	in a simplex facet will lie in $\mathcal{P}^d_{2s}$ and it satisfies
	$f(\pi(P)\cup Q)=f(P)\plusone f(Q)$.
\end{proposition}

\begin{corollary}
	The $f$-vector set of the simplicial $d$-polytopes forms an affine semigroup 
	\[
	\big(\fset(\mathcal{P}^d_{s}),\plusone\big),
	\]
	where the modified addition can be realized by 
	``glueing in simplex facets''.
\end{corollary}

We note that, just like $\pluszero$, the $\plusone$ operation does not
induce a semi-group structure $(\fset(\mathcal{P}^4),\plusone)$
on the $f$-vector set of all $4$-polytopes:
The data set of all $f$-vectors of $4$-polytopes with
$f_0+f_3\le22$
yields one single pair of $f$-vectors of $4$-polytopes
whose $\plusone$-sum is not an $f$-vector:

\begin{example}\label{example:non-additive2}
	Take $P=\Delta_4$ and
	$Q=(P^0_{2,2})^*=\Delta_2\times\Delta_2$,
	with $f$-vectors
	$f(P)=(5, 10, 10, 5)$ and
	$f(Q)=(9, 18, 15, 6)$.
	In this case there is no $4$-polytope with $f$-vector 
	$f(P)\plusone f(Q)=(10, 22, 21, 9)$.
\end{example}

We now turn to the connected sum operation.
Again, there is a large class of $d$-polytopes $\mathcal{P}^d_{cso}$,
those having a simple vertex $v$ and a simplex facet $F$ with $v\notin F$,
that is closed under the operation. In this case the class is also
closed under duality.

\begin{proposition}
	The $f$-vector set of the $3$-polytopes forms an affine semigroup 
	\[
	\big(\fset(\mathcal{P}^3),\plustwo\big),
	\]
	where the modified addition can be realized by 
	connected sum operations as follows:
	For any two $3$-polytopes $P,Q$ there are 
	$3$-polytopes $\bar{P},\bar{Q}$ with the same $f$-vectors
	such that the connected sum $\bar{P}\#\bar{Q}$ is a $3$-polytope with
	$f(\bar{P}\#\bar{Q})=f(P)\plustwo f(Q)$.
\end{proposition}

\begin{proof}
	Every $f$-vector of a $3$-dimensional polytope $P$ is also
	the $f$-vector of a $3$-dimensional polytope $P'$ with two     
	triangle facets: Indeed, every pyramid over an $n$-gon
	$P_n=v*C_n$ has two triangle facets, and this property
	is preserved under ``stacking onto a triangle facet''
	and under ``cutting off a simple vertex.''
	To prove Steinitz's lemma, see \cite{Stei3} or \cite[pp.~190/191]{Gr1-2}, 
	one shows that the resulting polytopes yield all the $f$-vectors of $3$-polytopes.
\end{proof}

We end this section with a conjecture, which we dare to state only
for $d=4$, also in view of Example~\ref{example:join_Cn*Cn} below.

\begin{conjecture}\label{conjecture:semigroup_d=4}
	The $f$-vector set of all $4$-polytopes is closed under the
	modified addition $\plustwo$. That is, 
	$\big(\fset(\mathcal{P}^4),\plustwo\big)$
	is an affine semigroup.
\end{conjecture}

There is a little bit of evidence for this:
The data set of all $184$ $f$-vectors of $4$-polytopes
with $f_0+f_3\le22$, which also produced Example~\ref{example:non-additive},
does not contain a counterexample to Conjecture~\ref{conjecture:semigroup_d=4}.
Also the $2$-dimensional coordinate projections of 
$\fset(\mathcal{P}^4)$, as determined by Grünbaum \cite[Sec.~10.4]{Gr1-2},
Barnette \& Reay \cite{barnette73:_projec}, and Barnette \cite{barnette74:_e_s},
do not produce any contradictions. (See Bayer \& Lee \cite{BaLee} and Höppner \& Ziegler \cite{Z59}
for summaries of this line of work.)

\section{The \emph{f}-vectors of 4-polytopes form an approximate monoid}%
\label{sec:approximate_monoids}%

The goal of this section will be to establish that for dimension $4$,
even though the modified sum of two $f$-vectors 
$f(P)\pluszero f(Q)$ is not in general the $f$-vector of a $4$-polytopes,
it is nevertheless ``close'' to one.

For this, we first prove that any $4$-polytope 
or its dual has a small facet. Then we deduce that the polytope
can be modified locally (without changing the $f$-vector much)
such that after the modification it has a simple vertex and a simplex facet.
This modification we perform on both the polytopes $P$ and $Q$ in question,
and then we take a connected sum that
approximately realizes the sum of the $f$-vectors.

\subsection{Small facets for 4-polytopes or their duals}

\begin{proposition}\label{prop:4poly_small_facet}
Let $P$ be a $4$-dimensional polytope with $f_0\le f_3$.
Then $P$ has a facet $F$ with at most
\[ 
	2 \frac{f_0}{(f_3)^{1/3}} + 2
\]
vertices.	
\end{proposition}

\begin{proof}
	Let $G$ be the vertex-facet incidence graph of $P$,
	which is a bipartite graph on $f_0+f_3$ vertices and $f_{03}$ edges.
	This bipartite graph $G$ is $K_{3,3}$-free, as any three distinct facets 
	of a $4$-polytope intersect in a face of dimension at most $1$, 
	which cannot have $3$ distinct vertices.
	
	From this, the K\H{o}v\'ari--S\'os--Tur\'an theorem \cite{KST}
	yields
	\[	 
	f_{03} \le z(f_0,f_3;3,3) \le 2^{1/3} (f_0-2) (f_3)^{2/3} + 2f_3
	\]
	and so the average number of vertices per facet, $f_{03}/f_3$, is bounded by
	\[	 
	\frac{f_{03}}{f_3} < 2 \frac{f_0}{f_3{}^{1/3}} + 2. 
	\] 
	Thus $P$ has a facet $F$ with $n\le2 f_0 (f_3)^{-1/3} + 2$
	vertices.
\end{proof}

\begin{corollary}
	Let $P$ be a $4$-polytope, then either $P$ or its dual $P^*$ has a facet
	$F$ with few vertices:
	\[
	   f_0(F) < \min_{0\le i\le3}\, 2\,f_i(P)^{2/3}+2.
	\] 
	Thus if $P$ is a ``large'' $4$-polytope (which means that it has many vertices,
	or equivalently many facets, or equivalently all components of the $f$-vector are large),
	then at least one of $P$ and $P^*$ has a ``small'' facet,
	whose number of vertices is small compared to \emph{all} the face numbers~$f_i(P)$.
\end{corollary}

\begin{proof}
	$P$ or $P^*$ satisfies $f_0\le f_3$. We may assume it's $P$, otherwise dualize.
	From Proposition~\ref{prop:4poly_small_facet} we then get that 
	$P$ has a facet $F$ with a small number of $n=f_0(F)$ vertices. Indeed,
	we know that  
	\[
	n < 2 \frac{f_0}{f_3{}^{1/3}} + 2 \le 2 f_3{}^{2/3}+2
	\]
	as well as 
	\[
	n < 2 \frac{f_0}{f_3{}^{1/3}} + 2 \le 2 f_0{}^{2/3}+2,
	\]
	using $f_0\le f_3$ in both estimates.
	Furthermore, we have $f_1\ge2f_0$ and $f_2\ge2f_3$ for $4$-polytopes,
	so the number of vertices of $F$ is small compared to
	$\min_i	f_i(P)$.
\end{proof}

\subsection{Modification and glueing}

The existence of a small facet (for $P$ or $P^*$) implies that
we can ``locally modify'' $P$ in such a way that the modified
polytope $P_\vartriangle$ differs from $P$ only a little in the sense
that only few faces are affected, and the $f$-vector change is also small.

\begin{proposition}\label{prop:simplex_facet}
	Let $P$ be a $4$-polytope,
	then there is an adapter $4$-polytope $P_\vartriangle$ 
	such that the $f$-vectors of $P$ and of $P_\vartriangle$ differ only by a little\textup:
	\[
	|f_i(P_\vartriangle) - f_i(P)| 
	\le 6\,f_i(P)^{2/3}+16
	\qquad\textrm{for all }i,\ 0\le i\le 3.
	\]
\end{proposition}
	
\begin{proof}
	Using duality, we may assume that $P$ has a small facet 
	(with a small number of $n$ vertices,
	as specified by Prop.\ \ref{prop:4poly_small_facet}).
	
	Now perform on $P$ a stellar subdivision on $F$, this will component-wise
	add at most 
$(1,n,3n-6,2n-5)$ to the $f$-vector of $F$.
	Then perform a stellar subdivision on the smallest resulting facet,
	which in the worst case is a pyramid over a pentagon, as the
	facet is a $3$-polytope, which has a $2$-face with at most $5$ vertices.
	This adds at most another 
$(1,6,10,5)$ to the $f$-vector.
	This creates several (at least $4$) tetrahedron facets.
	If we now do a stellar subdivision on one of these, then
	we are sure to also have a simple vertex; this last step adds 
$(1,4,6,3)$
	to the $f$-vector.
	Thus the $f$-vector difference
	between $P$ and the resulting polytope $P_\vartriangle$ is at most $(3,n+10,3n+10,2n)$.
	Hence 
    $|f_i(P_\vartriangle) - f_i(P)|
\le 3n+10
 < 3(2\,f_i(P)^{2/3}+2)+10
 =  6\,f_i(P)^{2/3}+16$.
\end{proof}

\begin{theorem}\label{thm:approximate_semigroup_d=4}
	Let $P',P''$ be $4$-polytopes,
	then there is a $4$-polytope $Q:=P'_\vartriangle\#P''_\vartriangle$ such that
	the $f$-vector of $Q$ differs from the sum of the $f$-vectors of $P'$ and of $P''$
	only by a little:
	\[
	|f_i(Q) - (f_i(P')+f_i(P''))| 
	< 
	12\,f_i(Q)^{2/3}+33
	\qquad\textrm{for all }i,\ 0\le i\le 3.
	\]
\end{theorem} 

\begin{proof}
	From Proposition~\ref{prop:simplex_facet} we know that we
	can modify $P'$ and $P''$ into adapter polytopes $P'_\vartriangle$ and $P''_\vartriangle$,
	such that
	$|f_i(P'_\vartriangle) - f_i(P')| \le 6\,f_i(P')^{2/3}+16$,
	and analogously for $P''$.
	
	Now we take the connected sum of these two polytopes, and get
	\[
	f(P'_\vartriangle  \# P'_\vartriangle) = 
	f(P'_\vartriangle) +  f(P'_\vartriangle) - (1,0,0,1).
	\]
	Thus we get
	\begin{eqnarray*}
	|f_i(Q) - (f_i(P')+f_i(P''))| 
	&\le
	&|f_i(Q) - f_i(P'_\vartriangle)+f_i(P''_\vartriangle)| +
	|f_i(P' _\vartriangle) - f_i(P' )| +
	|f_i(P''_\vartriangle) - f_i(P'')|
	\\
	&<&
	1 +
	(6f_i(P' )^{2/3}+16) + 
	(6f_i(P'')^{2/3}+16)
	\\
	&\le&12 f_i(Q)^{2/3}+33,
	\end{eqnarray*}
	where in the last step we use that 
	$f_i(P')\le f_i(P'_\vartriangle)\le f_i(Q)$ for all $i$, and analogously for~$P''$.%
\end{proof}

\subsection{Consequences}

For vectors $v\in\N^d$ (with positive integer entries),
let $\min(v)$ denote the size of the smallest entry, and
$\max(v)=|v|_\infty$ size of largest entry (i.e., the maximum norm);
any other norm would work as well, as long as we consider the dimension to be $d$ fixed.

Here is a suggestion for a definition that may be useful 
also in similar contexts:

\begin{definition}
	A set $S\subset \N^d$ is an \emph{approximate affine semigroup}
    if $v,w\in S$ implies that $v+w$ is close to an element $u\in S$,
	where \emph{close} is defined by $|u-(v+w)| \ll \min(u)$.
\end{definition}

\begin{corollary}
	The $f$-vectors of $4$-polytopes form an approximate affine semigroup
	$\fset(\mathcal{P}^4)\subset\N^4$.
\end{corollary}

\begin{definition}
Let $S\subset\N^d$.
A \emph{limit direction} is a unit vector $u\in S^d$  
such that for every $\varepsilon>0$ there are arbitrarily large 
vectors in $v\in S$ with $|\frac1{|v|}v -u|<\varepsilon$.	
\end{definition}

\begin{corollary}
	The set of limit directions of $\fset(\mathcal{P}^d)$
	is convex for $d\le4$.
\end{corollary}

\subsection{Dimension 5 and above?}

The following example shows that 
the approach that we employed for $4$-polytopes
cannot naively extended to higher dimensions:
It is plainly not true that every $5$-polytope
or its dual has a small facet:

\begin{example}\label{example:join_Cn*Cn}
	The join of $n$-gons 
	$P=C_n*C_n$ is a self-dual $5$-polytope
	with $f$-vector
	\[ 
	f(P)=(2n,n^2+2n,2n^2+2,n^2+2n,2n).
	\]
	
	All its facets are of type $I*C_n$
	(that is, a two-fold pyramid over an $n$-gon) with $f$-vector
	\[
	f(F)=(n+2,3n+1,3n+1,n+2),
	\]
	so they contain more than half of the vertices of $P$,
	and thus they have about the same ``size'' as 
	$f_0(P)$ and~$f_4(P)$.
\end{example}

\subsubsection*{Acknowledgement}
Thanks to Margaret Bayer, Isabella Novik, Arnau Padrol, and Francisco Santos for helpful discussions
at MSRI, 
to Philip Brinkmann, Moritz Firsching, and Hannah Sjöberg for valuable comments from Berlin,
and additionally to Moritz Firsching for the computations that yielded Examples~\ref{example:non-additive}
and~\ref{example:non-additive2} and to Johanna Steinmeyer for the \emph{Tikz}-pictures.


\begin{thebibliography}{10}

	\bibitem{altshuler84}
	{\sc A.~Altshuler and L.~Steinberg}, {\em Enumeration of the quasisimplicial
	  $3$-spheres and $4$-polytopes with eight vertices}, Pacific J. Math., 113
	  (1984), pp.~269--288.

	\bibitem{altshuler85}
	\leavevmode\vrule height 2pt depth -1.6pt width 23pt, {\em The complete
	  enumeration of the $4$-polytopes and $3$-spheres with eight vertices},
	  Pacific J. Math., 117 (1985), pp.~1--16.

	\bibitem{BabsonChan3}
	{\sc E.~K. Babson and C.~Chan}, {\em Counting faces for cubical spheres modulo
	  two}, Discrete Math., 212 (2000), pp.~169--183.

	\bibitem{barnette74:_e_s}
	{\sc D.~W. Barnette}, {\em The projection of the $f$-vectors of $4$-polytopes
	  onto the {$(E,S)$}-plane}, Discrete Math., 10 (1974), pp.~201--216.

	\bibitem{barnette73:_projec}
	{\sc D.~W. Barnette and J.~R. Reay}, {\em Projections of $f$-vectors of
	  four-polytopes}, J. Combinatorial Theory, Ser.~A, 15 (1973), pp.~200--209.

	\bibitem{BaBi}
	{\sc M.~M. Bayer and L.~J. Billera}, {\em Generalized {D}ehn--{S}ommerville
	  relations for polytopes, spheres and {E}ulerian partially ordered sets},
	  Inventiones Math., 79 (1985), pp.~143--157.

	\bibitem{BaLee}
	{\sc M.~M. Bayer and C.~W. Lee}, {\em Combinatorial aspects of convex
	  polytopes}, in Handbook of Convex Geometry, P.~Gruber and J.~Wills, eds.,
	  North-Holland, Amsterdam, 1993, pp.~485--534.

	\bibitem{Bjo2}
	{\sc A.~Bj\"orner}, {\em Face numbers of complexes and polytopes}, in
	  Proceedings of the International Congress of Mathematicians (Berkeley CA,
	  1986), 1986, pp.~1408--1418.

	\bibitem{Bjo6}
	\leavevmode\vrule height 2pt depth -1.6pt width 23pt, {\em Partial unimodality
	  for $f$-vectors of simplicial polytopes and spheres}, in Jerusalem
	  Combinatorics '93, H.~Barcelo and G.~Kalai, eds., vol.~178 of Contemporary
	  Math., Providence RI, 1994, Amer. Math. Soc., pp.~45--54.

	\bibitem{BjornerLinusson}
	{\sc A.~Bj\"orner and S.~Linusson}, {\em The number of $k$-faces of a simple
	  $d$-polytope}, Discrete Comput. Geometry, 21 (1999), pp.~1--16.

	\bibitem{BlBl1}
	{\sc G.~Blind and R.~Blind}, {\em Gaps in the numbers of vertices of cubical
	  polytopes {I}}, Discrete Comput. Geometry, 11 (1994), pp.~351--356.

	\bibitem{Z157}
	{\sc P.~Brinkmann and G.~M. Ziegler}, {\em Small {\itshape f}-vectors of
	  $3$-spheres and of $4$-polytopes}.
	\newblock Preprint, October 2016, 19 pages,
	  \href{http://arxiv.org/abs/1610.01028}{arXiv:1610.01028}; {Mathematics of
	  Computation}, to appear.

	\bibitem{eckhoff:_combin}
	{\sc J.~Eckhoff}, {\em Combinatorial properties of $f$-vectors of convex
	  polytopes}.
	\newblock Unpublished manuscript, Dortmund 1985.

	\bibitem{eckhoff:_combin2}
	\leavevmode\vrule height 2pt depth -1.6pt width 23pt, {\em Combinatorial
	  properties of $f$-vectors of convex polytopes}, Normat, 54 (2006),
	  pp.~146--159.

	\bibitem{Gr1-2}
	{\sc B.~Gr{\"u}nbaum}, {\em Convex {P}olytopes}, vol.~221 of Graduate Texts in
	  Math., Springer-Verlag, New York, 2003.
	\newblock Second edition prepared by V. Kaibel, V. Klee and G. M. Ziegler
	  (original edition: Interscience, London 1967).

	\bibitem{Z49-3}
	{\sc M.~Henk, J.~Richter-Gebert, and G.~M. Ziegler}, {\em Basic properties of
	  convex polytopes}, in Handbook of Discrete and Computational Geometry, J.~E.
	  Goodman, J.~O'Rourke, and C.~Toth, eds., Chapman \& Hall/CRC Press, Boca
	  Raton, third~ed., 2017, ch.~15, pp.~383--413.
	\newblock In print.

	\bibitem{Hoehn-DissETH1953}
	{\sc W.~H\"ohn}, {\em {Winkel und Winkelsumme im n-dimensionalen euklidischen
	  Simplex}}, PhD thesis, ETH Zurich, 1953.
	\newblock \url{https://doi.org/10.3929/ethz-a-000089763}.

	\bibitem{Z59}
	{\sc A.~H\"oppner and G.~M. Ziegler}, {\em A census of flag-vectors of
	  $4$-polytopes}, in Polytopes -- Combinatorics and Computation, G.~Kalai and
	  G.~M. Ziegler, eds., vol.~29 of DMV Seminars, Birkh\"auser-Verlag, Basel,
	  2000, pp.~105--110.

	\bibitem{kalai87:_rigid_i}
	{\sc G.~Kalai}, {\em Rigidity and the lower bound theorem, {I}}, Inventiones
	  Math., 88 (1987), pp.~125--151.

	\bibitem{KST}
	{\sc T.~K\H{o}v\'ari, V.~S\'os, and P.~Tur\'an}, {\em On a problem of {K}.
	  {Z}arankiewicz}, Colloq. Math., 3 (1954), pp.~50--57.

	\bibitem{Klee3}
	{\sc V.~Klee}, {\em A combinatorial analogue of {P}oincar\'e's duality
	  theorem}, Canadian J. Math., 16 (1964), pp.~517--531.

	\bibitem{miller04:_combin_commut_algeb}
	{\sc E.~Miller and B.~Sturmfels}, {\em Combinatorial Commutative Algebra},
	  vol.~227 of Graduate Texts in Math., Springer-Verlag, New York, 2004.

	\bibitem{Z163}
	{\sc H.~Sj\"oberg and G.~M. Ziegler}, {\em Semi-algebraic sets of $f$-vectors}.
	\newblock Preprint, 11 pages, September 2017.

	\bibitem{Stei3}
	{\sc E.~Steinitz}, {\em {\"U}ber die {E}ulerschen {P}olyederrelationen}, Archiv
	  der Mathematik und Physik, 11 (1906), pp.~86--88.

	\bibitem{Werner-diss}
	{\sc A.~Werner}, {\em Linear Constraints on Face Numbers of Polytopes}, PhD
	  thesis, TU Berlin, 2009.
	\newblock iv+185 pages,
	  \href{http://opus.kobv.de/tuberlin/volltexte/2009/2263/}{\url{opus.kobv.de/tuberlin/volltexte/2009/2263/}}.

	\bibitem{Z35}
	{\sc G.~M. Ziegler}, {\em Lectures on {P}olytopes}, vol.~152 of Graduate Texts
	  in Mathematics, Springer-Verlag, New York, 1995.
	\newblock Revised edition, 1998; seventh updated printing 2007.

	\bibitem{Z82}
	\leavevmode\vrule height 2pt depth -1.6pt width 23pt, {\em Face numbers of
	  $4$-polytopes and $3$-spheres}, in Proceedings of the International Congress
	  of Mathematicians (ICM 2002, Beijing), L.~Tatsien, ed., vol.~III, Beijing,
	  China, 2002, Higher Education Press, pp.~625--634.

	\end{thebibliography}
	\providecommand{\noopsort}[1]{}\providecommand{\noopsort}[1]{}

\end{document}